\documentclass[12pt]{article}
\usepackage{amsmath}
\usepackage{amssymb}

\newtheorem{theo}{Theorem}[section]
\newtheorem{lm}{Lemma}[section]

\newtheorem{rmk}{Remark}[section]

\allowdisplaybreaks

\numberwithin{equation}{section}

\def\Z{{\mathbb Z}}

\def\C{{\mathfrak C}}
\def\E{{\mathbf E}}
\def\Pg{{\mathbf P}}
\def\Po{{\mathbb P}}
\def\Prw{{\mathrm P}_\omega}
\def\1{{\mathbf 1}}
\def\Var{\mathop{\rm Var}\nolimits}

\def\sgn{\mathop{\rm sgn}}
\def\eps{\varepsilon}
\let\phi=\varphi
\def\qed{\hfill\rule{.2cm}{.2cm}}

\begin{document}

\title{Random walk  attracted by percolation clusters}

\author{S.~Popov\thanks{Partially supported by CNPq (302981/2002--0)}$^{~,1}$ \and
 M.~Vachkovskaia\thanks{Partially supported by CNPq (306029/2003--0)}$^{~,2}$ }

\maketitle

{\footnotesize
\noindent
$^1$ Departamento de Estat\'\i stica, Instituto de Matem\'atica e
Estat\'\i stica, Universidade de S\~ao Paulo, rua do Mat\~ao 1010,
CEP 05508--090, S\~ao Paulo SP, Brasil\\
e-mail: popov@ime.usp.br

\noindent
$^2$ Departamento de Estat\'\i stica, Instituto de Matem\'atica,
Estat\'\i stica e Computa\c{c}\~ao Cien\-t\'\i\-{}fica,
Universidade de Campinas,
Caixa Postal 6065, CEP 13083--970, Campinas SP, Brasil\\
e-mail: marinav@ime.unicamp.br

}

\begin{abstract}
Starting with a percolation model in $\Z^d$ in the subcritical regime,
we consider a random walk described as follows: the probability of
transition from~$x$ to~$y$ is proportional to some function~$f$ of
the size of the cluster of~$y$. This function is supposed to be increasing,
so that the random walk is attracted by bigger clusters. For $f(t)=e^{\beta t}$
we prove that there is a phase transition in~$\beta$, i.e., the random walk
is subdiffusive for large~$\beta$ and is diffusive for small~$\beta$.
\\[0.3cm]
\end{abstract}

\section{Introduction and results}
First, we describe the usual site percolation model in $\Z^d$.
This model is defined as follows.
For fixed $p\in (0,1)$, consider i.i.d.\ random variables $\omega(x)$, $x\in \Z^d$,
where $\omega(x)=1$ with probability~$p$ and $\omega(x)=0$ with probability
$1-p$. A site~$x$ is said to be open if $\omega(x)=1$ and closed otherwise.
Write $x\sim z$ if~$x$ and~$z$ are neighbors.
A (self-avoiding) path from~$x$ to~$y$ is: $\gamma(x, y)=\{x_0=x, x_1, x_2, \ldots, x_n=y\}$,
where $x_i\ne x_j$ if $i\ne j$ and $x_i\sim x_{i+1}$, $i=0, \ldots, n-1$.
A path~$\gamma$ is said to be open if all the sites in~$\gamma$ are open.
The cluster of~$x$ is defined by
\[
\C(x)=\{y\in \Z^d: \omega(y)=1 \mbox{ and there is an open path $\gamma(x,y)$
from $x$ to $y$}\}.
\]
Note that, if $\omega(x)=0$, then $\C(x)=\emptyset$.
It is a well-known fact (see e.g.~\cite{Gr}) that there exists
$p_{cr}$ (depending on~$d$; obviously, $p_{cr}=1$ in dimension 1) such that
if $p<p_{cr}$, then a.s.\ there is no infinite open cluster,
and if $p>p_{cr}$, then a.s.\ there exists an infinite open cluster
(also, with positive probability $|\C(0)|=\infty$).

Throughout this paper we assume that $p<p_{cr}$, i.e, the model is in the
(strictly) subcritical regime. Fix a parameter $\beta>0$.
The percolation configuration is regarded as random environment.
Fixed the environment, we start a discrete time random walk on $\Z^d$
with transition probabilities
\[
P^{\omega}_{xy}=\frac{e^{\beta |\mathfrak C(y)|}}{\sum\limits_{z\sim x} e^{\beta |\mathfrak C(z)|}},
\]
if $x\sim y$. Since~$\beta$ is positive, one can note that the random walk
is in some sense ``attracted'' by bigger clusters, and the strength of this attraction
grows with~$\beta$.
Denote by $\xi(t)$ the position of this random walk at time~$t$.
Let $\Po$ be the probability measure with respect to $\omega$ and $\Prw^x$ the (so-called
quenched) probability for the random walk starting from $x$ in the fixed
environment $\omega$. Denote also
$\Pg^x=\Prw^x \otimes \Po$
($\Pg^x$ is usually called the annealed probability); throughout the paper $\|\cdot\|$
stands for the $L_\infty$ norm.
Our main result is that there is a phase transition in $\beta$, i.e., the random walk
exhibits different behaviors for large and small $\beta$:
it is diffusive for small values of $\beta$
and subdiffusive for large values of $\beta$.
\begin{theo}
\label{main}
Suppose that the random walk $\xi(t)$ starts from the origin.
There exist $\beta_0$ and $\beta_1$ (depending on $d$) such that
$0<\beta_0\leq \beta_1 <\infty$ and
\begin{itemize}
\item[(i)] if $0<\beta<\beta_0$, then
\begin{equation}
\label{beta_small}
\lim_{t\to\infty}\frac{\log \max_{0\le s\le t}\|\xi(s)\|}{\log t}=\frac{1}{2}, \quad \Pg^0\mbox{-a.s.}
\end{equation}
\item[(ii)] if $\beta>\beta_1$, then
\begin{equation}
\label{beta_large}
\limsup_{t\to\infty}\frac{\log \max_{0\le s\le t}\|\xi(s)\|}{\log t}<\frac{1}{2}, \quad \Pg^0 \mbox{-a.s.}
\end{equation}
\end{itemize}
\end{theo}

One can prove also that the same  result holds for the bond percolation model
in the subcritical regime. The method of the proof remains the same; the reason why we have
chosen the site percolation is that for bond percolation there are some technical difficulties
(easily manageable, though; they relate to the fact that, in the bond percolation model,
two neighboring sites can belong to different large clusters)
in the proof of the part (ii) of Theorem~\ref{main}.

Recently much work has been done on the (simple or not) random walk on the unique infinite cluster
for the supercritical (bond or site) percolation in $\Z^d$ (see e.g.~\cite{B,BGP,GKZ,S};
see also~\cite{K} for some results for the random walk on the incipient infinite cluster
in dimension 2).
Another related subject is the class of models (see e.g.~\cite{BD,FM}) that can be described
as follows. Into each edge of $\Z^d$ we place a random variable that represents
the transition rate between the sites. The new features of the model of the
present paper are, first, the fact that the random environment is not independent,
and secondly, the absence of the uniform ellipticity. Speaking of uniform ellipticity,
we should mention that
in the paper~\cite{FIN} there was considered a simple symmetric
one-dimensional random walk with random rates,
 where the time spent at site $i$ before taking a step
has an exponential distribution with mean $\tau_i$, and $\tau_i$'s are i.i.d.\
positive random variables with distribution function $F$ having a polinomial tail.
One may find that there are similarities of the $d$-dimensional
analog of the model of~\cite{FIN} with our model, because clusters
of size $n$ will have ``density'' $e^{-Cn}$, and the mean time spent there
is roughly $e^{\beta n}$, so, thinking of clusters as ``sites'',
we indeed obtain a polinomial tail
of mean time spent at a given site. However, the facts that the random environment
is no longer independent and that here the random walk is not a time-change
of the simple random walk make the model of the
present paper considerably more difficult to analyze.

\section{Proof of Theorem~\ref{main}}
We begin by introducing some notations and recalling a well-known fact
from the percolation theory.
Namely, we will use the following result
(see~\cite{Men, Gr}): if $p<p_{cr}$, then there exists $c_1>0$ such
that for all $N>0$ and $x\in \Z^d$
\begin{equation}
\label{exp_oc} \Po[|\mathfrak C(x)|>N]\le e^{-c_1 N}.
\end{equation}
Now, to prove Theorem~\ref{main}, an important idea is to consider
$\xi(t)$ in finite region. Take $\Lambda_n=(-n/2, n/2]^d$ and
let the process $\xi^{(n)}(t)$ be the random walk $\xi(t)$ restricted on $\Lambda_n$.

\medskip
\noindent
{\it Proof of part (i).}\/
It can be easily seen that
$\xi(t)$ is reversible with the reversible measure
\begin{equation}
\label{rev1}
\pi(x)=e^{\beta |\mathfrak C(x)|}\sum\limits_{z\sim x} e^{\beta |\mathfrak C(z)|},
\end{equation}
and thus the finite Markov chain $\xi^{(n)}(t)$ is also reversible,
with the invariant (and reversible) measure
\begin{equation}
\label{rev}
\pi^{(n)}(x)=\frac{e^{\beta |\mathfrak C(x)|}\sum\limits_{z\sim x} e^{\beta |\mathfrak C(z)|}}{Z},
\end{equation}
where
\[
Z=\sum_{x\in\Lambda_n} e^{\beta |\mathfrak C(x)|}\sum_{z\sim x} e^{\beta |\mathfrak C(z)|}
\]
is the normalizing  constant, so that $\sum_{x\in\Lambda_n} \pi^{(n)}(x)=1$.

Consider also a random walk $\hat \xi^{(n)}(t)$ that is a continuization
of  $\xi^{(n)}(t)$. That is,  $\hat \xi^{(n)}(t)= \xi^{(n)}(N_t)$,
where $N_t$ is a Poisson process with rate $1$, independent of anything else
(in other words, $\hat \xi^{(n)}$ is a continuous time Markov chain with
the transition rates
equal to the transition probabilities of $\xi^{(n)}$).
Let ${\cal T}_i$ be the time interval between the jumps~$(i-1)$ and~$i$ of~$N_t$,
$S_n=\sum_{i=1}^n {\cal T}_i$, and $\hat T_A$ (respectively, $T_A$) be hitting time
of set~$A$ by random walk  $\hat \xi^{(n)}(t)$ (respectively, $\xi^{(n)}(t)$).
It can be easily seen (cf., for example, Chapter 2 of~\cite{AF})
that $\hat T_A=S_{T_A}$ and $\E(\hat T_A\mid T_A)=T_A$. Moreover, since $t^{-1}N_t \to 1$
a.s., many other results concerning $\hat \xi^{(n)}$ can be easily translated into the
corresponding results for $\xi^{(n)}$.

\begin{rmk}
\label{contin}
Using this technique, it is elementary to obtain that
Theorem~\ref{main} holds for~$\xi(t)$ iff it holds for~$\hat \xi(t)$,
where~$\hat \xi(t)$ is the continuization of~$\xi(t)$ defined in the same way.
\end{rmk}

So, now we consider the finite continuous time Markov chain $\hat \xi^{(n)}(t)$.
Denote by~$\lambda$ the spectral gap of $\hat \xi^{(n)}(t)$.
\begin{lm}
\label{lema_gap}
There exist $c_{14}>0$ and $n^*=n^*(\omega)$ such that for all $n>n^*$
we have $\lambda \ge c_{14} n^{-2}$, $\Po$-a.s.
\end{lm}

\medskip
\noindent
{\it Proof of Lemma~\ref{lema_gap}}.
The idea is to use Theorem~3.2.1 from~\cite{SC} to prove the lemma.
For each pair $x,y\in\Lambda_n$, we will   choose exactly one path $\gamma(x,y)$
(connecting~$x$ and~$y$) in a way described below.
Let $|\gamma(x,y)|$ be the length
of $\gamma (x,y)$ (i.e.\ the number of edges in $\gamma(x,y)$).
Denote by ${\cal E}(\Lambda_n)$ the set of edges of $\Z^d\cap \Lambda_n$.
For an  edge
$u=\langle z_1, z_2 \rangle$ let
$Q(u)=(P^{\omega}_{z_1 z_2}\pi^{(n)}(z_1)+P^{\omega}_{z_2 z_1}\pi^{(n)}(z_2))/2$.
According to Theorem~3.2.1 of~\cite{SC}, it holds that $\lambda \ge 1/A$, where
\begin{equation}
\label{gap1}
A=\max_{u\in {\cal E}(\Lambda_n)}\Big\{
\frac{1}{Q(u)}\sum _{x,y\in\Lambda_n: \;\gamma(x,y)\ni u}|\gamma(x,y)|\pi^{(n)}(x)\pi^{(n)}(y) \Big\}.
\end{equation}
Here,   we have, for $u=\langle z_1, z_2\rangle$,
\[Q(u)=\frac{e^{\beta(|\mathfrak C(z_1)|+|\mathfrak C(z_2)|)}}{Z}\]
and for each pair $x=(x^{(1)}, \ldots, x^{(d)}), y=(y^{(1)}, \ldots, y^{(d)}) \in\Lambda_n$
we choose the path $\gamma(x,y)$ in the following way.
Let $e_1, \ldots, e_d$ be the coordinate vectors. Denote $\Delta_i=y^{(i)}-x^{(i)}$
and let $\sgn(\Delta_i)$ be the sign of $\Delta_i$.
Suppose for definiteness that  $x^{(d)}\le y^{(d)}$ (so that $\sgn(\Delta_d)\ge 0$).
We take then
\begin{eqnarray*}
\gamma(x,y)&=&(x, x+e_d, \ldots, x+\Delta_d e_d, x+\sgn(\Delta_{d-1})e_{d-1}+\Delta_d e_d, \ldots,\\
&&x+\Delta_{d-1} e_{d-1}+\Delta_d e_d, \ldots, x+\Delta_1 e_1 +\cdots+\Delta_d e_d=y),
\end{eqnarray*}
so, first we successively change the $d$-th coordinate of $x$ to obtain the
$d$-th coordinate of $y$, then we do the same with $(d-1)$-th coordinate, and so on.
With this construction it is clear that the length of $\gamma(x,y)$
is at most~$dn$.
For an edge
\[
u=\langle (x^{(1)},\ldots,  x^{(d)}), (x^{(1)}, \ldots, x^{(d-1)}, x^{(d)}+1)\rangle
\]
define
\[
I_u=\{(z^{(1)},\ldots,  z^{(d)})\in\Lambda_n: \;\;
    z^{(1)}=x^{(1)},\ldots, z^{(d-1)}=x^{(d-1)},\; z^{(d)}\le x^{(d)}\}
\]
and
\[
R_u=\{(z^{(1)},\ldots, z^{(d)})\in\Lambda_n: \;\; z^{(d)}> x^{(d)}\}
\]
(for the edges of other directions the computations are quite analogous).
We have then
\begin{eqnarray}
\lefteqn{\sum_{x,y\in\Lambda_n: \;\gamma(x,y)\ni u}|
         \gamma(x,y)|\pi^{(n)}(x)\pi^{(n)}(y)}\phantom{****}\nonumber\\
&\le&dn\sum_{x\in I_u}\pi^{(n)}(x)\sum_{y\in R_u}\pi^{(n)}(y)\nonumber\\
&\le& dn\sum_{x\in I_u}\pi^{(n)}(x) \label{sum1},
\end{eqnarray}
as $\sum_{y\in R_u}\pi^{(n)}(y)\le 1$.

Now our goal is to  prove that with large probability,  for all such $u$,
$\sum_{x\in I_u}\pi^{(n)}(x)$ is of order $n/Z$. Denote
\[
\tilde I_u=\{(z^{(1)}, \ldots, z^{(d)})\in\Lambda_n: \;\;  z^{(1)}=x^{(1)},\ldots, z^{(d-1)}=x^{(d-1)}\}.
\]
Note that $I_u\subset \tilde I_u$, so
we will concentrate on obtaining an upper bound for
 $\sum_{x\in \tilde I_u}\pi^{(n)}(x)$.
It is important to observe that the variables $\pi^{(n)}(x)$
 are not independent.
For the sake of simplicity, suppose that~$\sqrt{n}$ is an integer, the general
case can be treated analogously.
Divide~$\tilde I_u$ into~$\sqrt{n}$ equal (connected) parts  of size~$\sqrt{n}$,
denote $x_{ij}:=(x^{(1)}, \ldots, x^{(d-1)}, (i-1)\sqrt{n}+j )$ and write
\[
\sum_{x\in \tilde I_u}\pi^{(n)}(x)= \sum_{j=1}^{\sqrt{n}}
              \sum_{i=-\frac{\sqrt{n}}{2}+1}^{\frac{\sqrt{n}}{2}} \pi^{(n)}(x_{ij}).
\]
For fixed~$j$, let $B_i=\1_{\{|\mathfrak C(
x_{ij})|<\sqrt{n}/2\}}$, and $B=\bigcap_{i=1}^{\sqrt{n}} B_i$.
By~(\ref{exp_oc}),  we have $\Po[B_i]\ge 1-e^{-c_1 \sqrt{n}}$  which
implies that
 $\Po[B]\ge 1-e^{-c_2 \sqrt{n}}$ for some $c_2>0$,
so
\begin{eqnarray*}
&&\Po\Big[\sum_{i=-\frac{\sqrt{n}}{2}+1}^{\frac{\sqrt{n}}{2}} \pi^{(n)}(x_{ij})\ge \frac{k\sqrt{n}}{Z}\Big]\\
&&~~~~~\le
\Po\Big[\Big(\sum_{i=-\frac{\sqrt{n}}{2}+1}^{\frac{\sqrt{n}}{2}} \pi^{(n)}(x_{ij})\ge
                         \frac{k\sqrt{n}}{Z}\Big)\1_B\Big]+e^{-c_2 \sqrt{n}}.
\end{eqnarray*}
Now, it is important to note that the variables
$e^{\beta |\mathfrak C(x_{ij})|}\sum_{z\sim x_{ij}} e^{\beta |\mathfrak C(z)|}\1_B =
\pi(x_{ij})\1_B$
(recall~(\ref{rev1}) and~(\ref{rev})), $i=1, \ldots, \sqrt{n}$,
are independent. We have also
\[
e^{\beta |\mathfrak C(x_{ij})|}\sum_{z\sim x_{ij}} e^{\beta |\mathfrak C(z)|}\le
                                                    2d e^{2\beta |\mathfrak C(\tilde x_{ij})|},
\]
where $\tilde x_{ij}$ satisfies $|\mathfrak C(\tilde x_{ij})|=
                     \max_{z\sim x_{ij}} \{|\mathfrak C(x_{ij})|, |\mathfrak C(z)|\}$.

For  $y=n^\alpha$ ($\alpha>0$ will be chosen later), using~(\ref{exp_oc}), we have
\begin{equation}
\label{N1ocenka}
\sum_{i=-\frac{\sqrt{n}}{2}+1}^{\frac{\sqrt{n}}{2}}
                       \Po[2d e^{2\beta |\mathfrak C(\tilde x_{ij})|}>y]
                         \le \sqrt{n} e^{-\frac{c_3 \log y}{2\beta}}
=\sqrt{n} y^{-\frac{c_3}{2\beta}}= n^{\frac{1}{2}-\frac{c_3\alpha}{2\beta}}.
\end{equation}

According to Corollary~1.5 from~\cite{N},
if $X_1, \ldots, X_k$ are independent random variables, $S_k=\sum_{i=1}^k X_i$, $F_i(x)=\Po[X_i<x]$, then
for any set $y_1, \ldots, y_k$ of positive numbers and any~$t$, $t\in(0,1]$,
\begin{equation}
\label{N1}
\Po[S_k\ge x]\le \sum_{i=1}^k \Po[X_i>y_i]+\Big(\frac{eA_t^+}{xy^{t-1}}\Big)^{\frac{x}{y}},
\end{equation}
where
$y\ge\max\{y_1, \ldots, y_k\}$ and
\[A_t^+=\sum_{i=1}^k\int\limits_0^{\infty}u^t dF_i(u).\]

Denote $\bar F_i(u)=1-F_i(u)$. We apply Corollary~1.5 from~\cite{N}
to random variables $2d e^{2\beta |\mathfrak C(\tilde x_{ij})|} \1_B$,
$i=-\frac{\sqrt{n}}{2}+1, \ldots, \frac{\sqrt{n}}{2}$,  with
$x=k\sqrt{n}$, $y_i\equiv y=n^{\alpha}$,
and $t=1$. First term of the right-hand side of~(\ref{N1})
was estimated in~(\ref{N1ocenka}). For the second term, we have
\begin{eqnarray}
A_1^+&=&\sum_{i=-\frac{\sqrt{n}}{2}+1}^{\frac{\sqrt{n}}{2} }\int\limits_0^{\infty}u dF_i(u)\nonumber\\
&=&-\sqrt{n}\int\limits_0^{\infty}u d\bar F_i(u)\nonumber\\
&=&-\sqrt{n}\int\limits_0^{\infty} \bar F_i(u) du\nonumber\\
&\le &c_4\sqrt{n} \int\limits_1^{\infty} u^{-\frac{c_3}{2\beta}}du\nonumber\\
&=& \label{N2ocenkaA}c_5 \sqrt{n},
\end{eqnarray}
as $\beta$ is small, thus
\begin{equation}
\label{N2ocenka}
\Big(\frac{eA_t^+}{xy^{t-1}}\Big)^{\frac{x}{y}}\le
        \Big(\frac{e c_5\sqrt{n}}{k \sqrt{n}}\Big)^{\frac{k\sqrt{n}}{n^{\alpha}}}=
    \Big(\frac{e c_5 }{k}\Big)^{k n^{\frac{1}{2}-\alpha}}
\end{equation}
so, to guarantee that $ \Big(\frac{eA_t^+}{xy^{t-1}}\Big)^{\frac{x}{y}}\to 0$ as $n\to\infty$,
  it is sufficient to take
$k$ large enough and $\beta/c_3<\alpha< 1/2$.

We proved that for $\beta$ sufficiently small
\begin{equation}
\label{Nocenkafinal}
\Po\Big[\sum_{i=-\frac{\sqrt{n}}{2}+1}^{\frac{\sqrt{n}}{2}} \pi^{(n)}(x_{ij})\ge \frac{k\sqrt{n}}{Z}\Big]\le
                          c_7 n^{-\frac{1}{2}(\frac{c_3\alpha}{\beta}-1)}
\le c_7 n^{-\frac{c_{10}}{\beta}+\frac{1}{2}}.
\end{equation}
Thus,
\begin{eqnarray}
\Po\Big[\sum_{x\in \tilde I_u}\pi^{(n)}(x)\ge  \frac{c_8 n}{Z}\Big]
      &=&\Po\Big[\sum_{j=1}^{\sqrt{n}}\sum_{i=-\frac{\sqrt{n}}{2}+1}^{\frac{\sqrt{n}}{2}}
                                                      \pi^{(n)}(x_{ij})\ge  \frac{c_8n}{Z}\Big]\nonumber\\
&\le& \sum_{j=1}^{\sqrt{n}} \Po\Big[\sum_{i=-\frac{\sqrt{n}}{2}+1}^{\frac{\sqrt{n}}{2}}
                                                       \pi^{(n)}(x_{ij})\ge \frac{k\sqrt{n}}{Z}\Big]\nonumber\\
&\le& \label{sum_total_pi} c_7 n^{-\frac{c_{10}}{\beta}+1}
\end{eqnarray}

 So, using~(\ref{sum_total_pi}) in~(\ref{sum1}) and~(\ref{gap1}), we have
\begin{equation}
\label{gap2}
A\le \max_{u\in{\cal E}(\Lambda_n)}
 \Big\{\frac{c_9Z}{e^{\beta(|\mathfrak C(z_1)|+|\mathfrak C(z_2)|)}}  \frac{n^2}{Z}\Big\}\le c_9 n^2,
\end{equation}
as $e^{\beta(|\mathfrak C(z_1)|+|\mathfrak C(z_2)|)}\ge 1$,
with probability at least $1-c_7 n^{-\frac{c_{10}}{\beta}+1}$. Since  $\beta$
can be made arbitrarily small, Borel-Cantelli lemma implies that
for almost all environments for $n$ large enough it holds that $A\le c_{13} n^2$
and thus $\lambda\ge c_{14}n^{-2}$. Lemma~\ref{lema_gap} is proved.
\qed

\medskip

Now, using Lemma~2.1.4 from~\cite{SC} with $f(x)=\1_{\{\|x\|\ge n/4\}}$, where,
as before, $\|\cdot\|$ is the $L_{\infty}$ norm, we prove~(\ref{beta_small}).
By Lemma~2.1.4 from~\cite{SC}
we have that
\[
\|H_t f-\pi^{(n)}(f)\|_2^2\le e^{-2\lambda t} \Var_\pi^{(n)}(f).
\]
In what follows we
 show that $\pi^{(n)}(f)$ is of constant order.
We have
\[
\pi^{(n)}(f)=\sum_{x\in\Lambda_n, \; \|x\|\ge n/4}\frac{e^{\beta |\mathfrak C(x)|}\sum\limits_{z\sim x}
                 e^{\beta |\mathfrak C(z)|}}{Z}.
\]
Since $|\mathfrak C(x)|\ge 0$, it is easy to obtain that  for all
$\omega$ it holds
\begin{equation}
\label{sum_ring1}
\sum_{x\in\Lambda_n, \; \|x\|\ge n/4}e^{\beta |\mathfrak C(x)|}\sum\limits_{z\sim x}
                        e^{\beta |\mathfrak C(z)|}\ge \frac{n^d}{2}.
\end{equation}
Using the same kind of argument as in the proof of Lemma~\ref{lema_gap},
one can easily see that for all~$n$
\begin{equation}
\label{sum_ring2}
\Po\Big[\sum_{x\in\Lambda_n, \; \|x\|< n/4}e^{\beta |\mathfrak C(x)|}\sum\limits_{z\sim x}
             e^{\beta |\mathfrak C(z)|}\ge c_{15}n^d\Big]
     \le c_{15}'' n^{-\frac{c_{15}'}{\beta}},
\end{equation}
where $c_{15}',c_{15}''$ depend only on $c_{15}$.
Thus, with probability at least $1-c_{15}''n^{-\frac{c_{15}'}{\beta}}$
we have $\pi^{(n)}(f)\ge const$.
Then, using that $\Var_\pi^{(n)}(f)\le 1$, taking $t=c_{16} n^2$
for~$c_{16}$ large enough
yields that the random walk $\hat \xi^{(n)}(t)$, and thus $ \xi^{(n)}(t)$,
will be at distance  of order $n$ from the origin (as both random walks start from $0$) after
a time of order $n^2$ with probability bounded away from $0$.

Now, for any fixed $\eps>0$, divide the time interval $(0,t]$ into $t^\eps$
intervals of length $t^{1-\eps}$. Borel-Cantelli lemma implies then that
for $t$ large enough there will be at least one time interval
such that at the end of this interval $\xi^{(n)}$ will be at distance at least
$t^{\frac{1}{2}-\frac{\eps}{2}}$ from the origin. Since $\eps>0$ is arbitrary,
we proved that
\begin{equation}
\label{i1}
\liminf_{t\to\infty}\frac{\log \max_{0\le s\le t}\|\xi(s)\|}{\log t}\ge\frac{1}{2},
           \quad \Pg^0 \mbox{-a.s.}
\end{equation}

It remains to prove that
\begin{equation}
\label{i2}
\limsup_{t\to\infty}\frac{\log \max_{0\le s\le t}\|\xi(s)\|}{\log t}\le\frac{1}{2},
          \quad \Pg^0 \mbox{-a.s.}
\end{equation}
It is a well-known fact that a reversible Markov chain with a ``well-behaved"
reversible measure cannot go much farther than $t^{1/2}$ by time $t$,
see~\cite{BP, C, K, V}.
By Theorem~1 from~\cite{C}
we have, for any $\eps>0$
\begin{eqnarray}
\Prw[\|\xi_n\|\ge n^{1/2+\eps}]&\le &2e^{-\frac{n^{2(1/2+\eps)}}{2n}}\sum_{y: \|y\|\ge n^{1/2+\eps}}
   \frac{e^{\beta|\C(y)|}\sum_{z'\sim y}e^{\beta|\C(z')|}}{ e^{\beta|\C(0)|}\sum_{z\sim 0}e^{\beta|\C(z)|}}
       \nonumber\\
&\le & c_{20} e^{-n^{\eps}} \quad \Po\mbox{ -a.s.} \label{oc_sverhu_i}
\end{eqnarray}
for all $n$ large enough. To obtain the bound~(\ref{oc_sverhu_i})
we have used the fact that, due to~(\ref{exp_oc}),
\[\Po[\max_{x\in \Lambda_n}e^{\beta |\C(x)|}\ge n^{\frac{\eps}{2}}]\le
n^{-\frac{c_{21}\eps}{\beta}}\]
for some $c_{21}>0$.
Borel-Cantelli lemma and~(\ref{oc_sverhu_i}) imply~(\ref{i2})
and thus the  part (i) of Theorem~\ref{main} is proved.

\medskip
\noindent
{\it Proof of part~(ii)}. For $x\in \Z^d$
let $T_0(x)=0$, $T_0'(x)=0$, and define
\begin{eqnarray*}
T_i'(x)&=&\min\{t\ge T_{i-1}(x)+T_{i-1}'(x): \; \xi(t)\in \C(x)\}\\
T_i(x)&=&\min\{t>T_i'(x):\; \xi(t)\notin \C(x)\}-T_i'(x),
\end{eqnarray*}
$i=1,2,3, \ldots$, where $T_i(x)$ is defined if
$\min\{t\ge T_{k}(x): \; \xi(t)\in \C(x)\}$, $k\le i$, are finite.
In words, $T_i'(x)$ is the moment of $i$th entry to the cluster of~$x$,
and $T_i(x)$ is the time spent there (i.e., after $T_i'(x)$ and
before going out of~$\C(x)$). It is important to
note that the cluster~$\C(x)$ is surrounded by sites with
$\omega(\cdot)=0$. Comparing~$T_i(x)$
with geometric random variable with parameter
$(2d-1+e^{\beta |\C(x)|})^{-1}$, one can easily see that if $|\C(x)|\ge \delta \log n$,
then $\E T_i(x)\ge c_{17}n^{\beta\delta}$. Moreover, it is elementary to
obtain that for $\eps>0$ and for any $\delta>0$ we can choose~$\beta$
large enough so that with probability bounded away from~$0$ we have
\begin{equation}
\label{zalipaet}
T_i(x)\ge  c_{18}n^{2+\eps}
\end{equation}
for $x$ such that $|\C(x)|\ge \delta \log n$.

Now, we use a dynamic construction of the percolation
environment usually called the {\it generation method\/} (see~\cite{Men}).
That is, we proceed in the following way:
we assign generation index~$0$ to the origin, and
put~$\omega(0)=1$ (the origin is open) with probability~$p$
or $\omega(0)=0$ (closed) with probability $1-p$.
If $\omega(0)=0$, the process stops.
If $\omega(0)=1$, then to all $x\sim 0$ we assign
generation index~$1$, and put $\omega(x)=1$  with probability~$p$
or $\omega(x)=0$  with probability $1-p$, independently.
Suppose that the we constructed~$m$ generations of the process. Let~$Y_i$
be the set of sites with generation index~$i$ and $Y^m=\{0\}\cup Y_1\cup\ldots\cup Y_m$.
Denote by~$Y_{m+1}$ the set of neighbors of the open sites in~$Y_m$ which do not belong
to~$Y^m$. Assign to the sites from $Y_{m+1}$ the generation index~$m+1$ and
a value~$1$ or~$0$ in a way described above. If $Y_m\ne \emptyset$ and $\omega(y)=0$
for all $y\in Y_{m+1}$, then the process stops.
Note that for subcritical percolation this process stops a.s., and what we obtain
at the moment when the process stops is the cluster of the origin surrounded by $0$-s.

So, first we construct the environment within the set $H_1=\C(0)\cup \partial \C(0)$, where
\[
\partial \C(0)=\{y: \; y\notin \C(0), y\sim x\mbox{ for some }x\in\C(0)\}
\]
(note that $\omega(y)=0$ for any $y\in \partial \C(0)$)
and we know nothing yet about the environment out of the set~$H_1$.
For an arbitrary set $H\subset\Z^d$ denote
\begin{eqnarray*}
H^\circ&=&\{x\notin H: \; \mbox{ for any infinite path
$\gamma(x)$ starting from $x$ it holds}\\
&&~~~~~~\mbox{that }\gamma(x)\cap H\ne\emptyset\}
\end{eqnarray*}
(i.e., $H^\circ$ is the set of the holes within the set $H$) and let
\[
G_1=H_1\cup H_1^\circ .
\]
Then, choose $\omega (x)$  for $x\in H_1^\circ$ and start the random walk $\xi(t)$
 from the origin. Let
\[
\tau_1=\min\{t:\; \xi(t)\notin G_1\}.
\]
Note that $\C(\xi(\tau_1))\cap \C(0)=\emptyset$,
and construct, using the above method
\[
H_2=\C(\xi(\tau_1))\cup \partial \C(\xi(\tau_1))
\]
and
\[
G_2=G_1\cup H_2\cup (G_1\cup H_2 )^\circ.
\]
Then, define
\[
\tau_2=\min\{t:\; \xi(t)\notin G_2\},
\]
and so on. For all $i$, we have
\[
 \Po[\C(\xi(\tau_i)) \ge \delta\log n]\ge p^{\delta\log n}
\]
where~$p$ is the percolation parameter.
This is so due to the fact that, to have $\C(\xi(\tau_i)) \ge \delta\log n$,
it is sufficient to choose a path of length $\delta\log n$
emanating from $\xi(\tau_i)$
which does not intersect~$G_i$ (it is possible by the
construction of $\tau_i$, since $\xi(\tau_i)$ cannot be completely
surrounded by points of~$G_i$),
and such path will be open with probability $p^{\delta\log n}$.
Thus, for any $\eps>0$ (one can take the same~$\eps$ from~(\ref{zalipaet})),
we can choose~$\delta$ small enough (take~$\delta$ such that $\delta\log p^{-1}<\eps$) so that
\begin{equation}
\label{lovitsya}
\Pg^0[\C(\xi(\tau_i)) \ge \delta\log n\mid {\cal F}_i]\ge n^{-\eps},
\end{equation}
where ${\cal F}_i$ is the $\sigma$-algebra generated by
$\{\omega(x), \; x\in G_i\}$ and $\{\xi(m), \; m\le \tau_{i}\}$.
Fix $\theta>0$ in such a way that $1-\theta>\eps$.
Note that, as $p<p_{cr}$, using~(\ref{exp_oc}) and Borel-Cantelli lemma,
for~$n$ large enough
$\min\{k:\; \xi(\tau_k)\notin \Lambda_n\}$ (the number of times that we
repeat the basic step in the above construction)
 will be of order at least $n^{1-\theta}$ for all~$n$ large enough,
$\Pg^0$-a.s.\ (recall that $\Lambda_n=(-n/2, n/2]^d$; with overwhelming probability
all the clusters inside $\Lambda_n$ will be of sizes at most~$n^\theta$).
On each step, by~(\ref{lovitsya}), with probability at least~$n^{-\eps}$ the random
walk enters the cluster of size at least~$\delta\log n$.
By~(\ref{zalipaet}), it stays in that cluster
(if~$\beta$ is large enough) for at least $n^{2+\eps}$ time units with large
probability. If $1-\theta>\eps$, with overwhelming probability
on some step (of the above construction) the random walk will delay
(in the corresponding cluster) for more
than $n^{2+\eps}$ time units before going out of~$\Lambda_n$.
In other words, we will have
\[
 \max_{s\leq n^{2+\eps}} \|\xi(s)\| \leq \frac{dn}{2},
\]
which implies~(\ref{beta_large}).
 This concludes the proof of Theorem~\ref{main}.
\qed


\begin{thebibliography}{19}
\bibitem{AF} {\sc D.~Aldous, J.A.~Fill}
{\it Reversible Markov Chains and Random Walks
on Graphs.} Available at: {\tt http://www.stat.berkeley.edu/users/} {\tt aldous/RWG/book.html}

\bibitem{B} {\sc M.T.~Barlow} (2004)
Random walk on supercritical percolation clusters.
{\it Ann. Probab.} {\bf 32} (4), 3024--3084.

\bibitem{BP} {\sc M.T.~Barlow, E.A.~Perkins} (1989)
Symmetric Markov chains in $\Z^d$: how fast can they move?
{\it Probab. Theory Relat. Fields} {\bf 82}, 95--108.

\bibitem{BGP} {\sc N.~Berger, N.~Gantert, Y.~Peres} (2003)
The speed of biased random walk on percolation clusters.
{\it Probab. Theory Relat. Fields} {\bf 126} (2), 221--242.

\bibitem{BD} {\sc D.~Boivin, J.~Depauw} (2003)
Spectral homogenization of reversible random walks on $\Z^d$ in a random environment.
{\it Stochastic Process. Appl.} {\bf 104}, 29--56.

\bibitem{C} {\sc T.K.~Carne} (1985)
A transmission formula for Markov chains.
{\it Bull. Sc. Math. (2)}, {\bf 109} (4), 399--405.

\bibitem{FIN} {\sc L.R.G.~Fontes, M.~Isopi, C.M.~Newman} (2002)
      Random walks with strongly inhomogeneous rates and singular
      diffusions: convergence, localization and aging in one dimension.
      {\it  Ann.~Probab.} {\bf 30}, 579--604.

\bibitem{FM} {\sc L.R.G.~Fontes, P.~Mathieu} (2005)
On symmetric random walks with random conductancies on $\Z^d$.
Preprint.

\bibitem{Gr} {\sc G.R.~Grimmett} (1999)
{\it Percolation.} Springer, Berlin.

\bibitem{GKZ} {\sc G.R.~Grimmett, H.~Kesten, Y.~Zhang} (1993)
Random walk on the infinite cluster of the percolation model.
{\it Probab. Theory Relat. Fields} {\bf 96} (1), 33--44.

\bibitem{K} {\sc H. Kesten} (1986)
Subdiffusive behaviour of random walk on a random cluster.
{\it Ann. Inst. Henri Poincar\'e} {\bf 22} (4), 425--487.

\bibitem{Men} {\sc M.V.~Menshikov} (1986)
 Coincidence of critical points
in percolation problems.  {\it Sov. Math. Doklady} {\bf  33}, 856--859.

\bibitem{N} {\sc S.V.~Nagaev} (1979)
Large deviations of sums of independent random variables.
{\it Ann. Probab.} {\bf 7} (5), 745--789.

\bibitem{SC} {\sc L.~Saloff-Coste} (1997)
{\it Lectures on Finite Markov Chains.}
Lectures on probability theory and statistics (Saint-Flour,
1996), 301--413, Lecture Notes in Math., {\bf 1665}, Springer, Berlin.

\bibitem{S} {\sc A.S.~Sznitman} (2003)
On the anisotropic random walk on the supercritical percolation cluster.
{\it Commun. Math. Phys.}, {\bf 240}, 123--148.

\bibitem{V} {\sc N.Th. Varopoulos} (1985)
Long range estimates for Markov chains.
{\it Bull. Sc. Math. (2)}, {\bf 109} (3), 225--252.

\end{thebibliography}
\end{document}